\documentclass[12pt]{article}
\usepackage{amsmath,amsfonts,amssymb}

\def\qed{\nopagebreak\hfill{\rule{4pt}{7pt}}}
\def\proof{\noindent {\it{Proof.} \hskip 2pt}}

\newtheorem{theo}{Theorem}[section]

\newtheorem{lemm}[theo]{Lemma}

\newtheorem{coro}[theo]{Corollary}

\begin{document}
\begin{center}
{\large \bf  Stanley's Zrank Problem on Skew Partitions}
\end{center}

\begin{center}
William Y. C. Chen$^{1}$ and Arthur
L. B. Yang$^{2}$\\[6pt]
Center for Combinatorics, LPMC\\
Nankai University, Tianjin 300071, P. R. China\\
Email: $^{1}${\tt chen@nankai.edu.cn}, $^{2}${\tt
yang@nankai.edu.cn}
\end{center}

\vspace{0.3cm} \noindent{\bf Abstract.} We present an affirmative
answer to Stanley's zrank problem, namely, the zrank and rank are
equal for any skew partition. We show that certain classes of
restricted Cauchy matrices are nonsingular and furthermore, the
signs depend  on the number of zero entries.   Similar to notion
of the jrank of a skew partition, we give a characterization of
the rank in terms of the Giambelli type matrices of the
corresponding skew Schur functions. We also show that the sign of
the determinant of a factorial Cauchy matrix is uniquely
determined by  the number of its zero entries, which implies the
nonsingularity of the inverse binomial coefficient matrix.



\section{Introduction}

In the study of tensor products of Yangian modules, Nazarov and
Tarasov \cite{NT} gave a generalization of the rank of an ordinary
 partition to a skew partition.   Stanley \cite{Stanley}
obtained several characterizations of the rank of a skew partition
in terms of the reduced partition code, the Jacobi-Trudi matrix,
and the minimal border strip decomposition. Stanley also
introduced the notion of the zrank of a skew partition in terms of
the specialization of the skew Schur function, and proposed the
problem whether the zrank and the rank are always equal.

Yan, Yang and Zhou \cite{YYZ} gave an equivalent characterization
of Stanley's problem in terms of the restricted Cauchy matrix
based on two integer sequences. In this paper, we extend the
definition of a restricted Cauchy matrix to two sequences of real
numbers subject to certain conditions. We prove that every
restricted Cauchy matrix is nonsingular, and thus give an
affirmative answer to Stanley's problem.

In the spirit of Stanley's notion of the jrank
 of a skew partition which is defined
 as the number of rows in the Jacobi-Trudi matrix
 in which one does not appear, we
 introduce the notion of {\em grank} in terms of the Giambelli type
 matrix defined by Hamel and Goulden
 for a skew Schur function \cite{HG}. Given any outside
decomposition of a skew partition,  the grank is defined by the
number of rows in which there are no entries equal to one. It
turns out that the grank is well-defined, namely, it does not
depend on the outside decomposition of the skew partition. We show
that the grank is always equal to the rank for any skew partition.

This paper is also concerned with the nonsingularity of the
factorial Cauchy matrices. Given a sequence $A$ of real numbers
and a sequence $B$ of integers, we define the factorial Cauchy
matrix to be a matrix with each entry being either the inverse of
the falling factorial or zero, similar to the definition of the
restricted Cauchy matrix. We prove that the determinant of the
factorial Cauchy matrix has the same property as the restricted
Cauchy matrix. A special case of the factorial Cauchy matrix falls
into the framework of the calculation of some determinants
involving the $s$-shifted factorial by Normand \cite{Normand} in
the study of the probability density of the determinant of random
matrices \cite{DL, MN}.

The double Schur functions serve as a tool for proving the
nonsingularity of the factorial Cauchy matrix without zero
entries. The double Schur functions are a natural extension of the
factorial Schur function introduced by Biedenharn and Louck
\cite{BL}, and further studied by Chen and Louck \cite{C-L},
Goulden and Greene \cite{GG}, Macdonald \cite{Macdonald},
 Chen, Li and Louck \cite{CLL}.

As a direct application of the nonsingularity of factorial Cauchy
matrices, we prove the nonsingularity of inverse binomial
coefficient matrices, which are defined as matrices with entries
being either zeros or the inverses of the binomial coefficients.

\section{The restricted Cauchy matrices}

Let $A=(a_1, \ldots, a_n)$ and $B=(b_1, \ldots, b_n)$ be two
sequences of real numbers. Suppose that $A$ is strictly
decreasing, $B$ is strictly increasing, and $a_i>b_{n+1-i}$ and
$a_i\neq b_j$ for any $i,j$. We define a matrix
${\textbf{C}}(A,B)=(c_{ij})_{i,j=1}^n$ by setting
\begin{equation}\label{def-eq} c_{ij}=\left\{
\begin{array}{ll}
{\displaystyle \frac{1}{a_i-b_j}}, & \mbox{ if $a_i>b_j$},\\[12pt]
0, & \mbox{ if $a_i<b_j$}.
\end{array}
\right.
\end{equation}

Similar to the definition in \cite{YYZ}, a matrix $M$ is called a
\emph{restricted real Cauchy matrix} if there exist two sequences
$A$ and $B$ satisfying the above conditions such that
$M=\textbf{C}(A,B)$. If $A$ and $B$ are restricted to integer
sequences, we call $M$ a \emph{restricted integer Cauchy matrix}.

Let $\omega(M)$ be the number of zero entries of $M$.  We have the
following criterion for the sign of the restricted Cauchy
determinant.

\begin{theo}\label{mainthm} Any restricted real Cauchy matrix
$M=\textbf{C}(A,B)$ is nonsingular. Furthermore, the determinant
$\det(M)$ is positive if $\omega(M)$ is even; or negative if
$\omega(M)$ is odd.
\end{theo}

Before giving the proof of the above theorem, let us  recall some
definitions on matrices.  A matrix $M=(m_{ij})_{i,j=1}^n$ is
called \emph{reducible} if the indices $1, 2, ..., n$ can be
divided into two disjoint nonempty subsets $\{i_1, i_2, ...,i_s\}$
and $\{j_1, j_2, ...,j_t\}$  (with $s+t=n$) such that
$$\mbox{$m_{i_{\alpha},j_{\beta}}=0$, for $\alpha=1, 2, ...,s$ and
$\beta=1, 2, ...,t$.}$$ Otherwise, $M$ is said to be an
irreducible matrix. Clearly, a restricted real Cauchy matrix
$M=\textbf{C}(A,B)$ is irreducible if $a_i>b_{n+2-i}$ for
$i=2,3,\ldots,n$. Given a square matrix $M=(m_{ij})_{i,j=1}^n$,
let $M_{ij}$ denote the $(i,j)$-th \emph{minor} of $M$ which is
the matrix obtained from $M$ by deleting the $i$-th row and the
$j$-th column. We have the following lemma:

\begin{lemm}\label{minorcauchy}
If $M=\textbf{C}(A,B)$ is an irreducible restricted Cauchy matrix,
then each minor $M_{ij}$ is a restricted Cauchy matrix.
\end{lemm}

\proof Let
\begin{eqnarray*}
A'&=&(a_1',\ldots,a_{n-1}')=(a_1,\ldots,\hat{a}_i,\ldots,a_n)\\
B'&=&(b_1',\ldots,b_{n-1}')=(b_1,\ldots,\hat{b}_j,\ldots,b_n)
\end{eqnarray*}
where $\hat{}$ stands for the  missing entry from the sequence. It
suffices to prove that $M_{ij}=\textbf{C}(A',B')$. Clearly, each
entry of $M_{ij}$ is defined by \eqref{def-eq}. We only need to
show that $a_{k}'>b_{n-k}'$ for any $k:1\leq k\leq n-1$. There are
four cases:
\begin{enumerate}
\item[(a)] If $k<i$ and $n-k<j$, then
$a_k'=a_k>b_{n+1-k}>b_{n-k}=b_{n-k}'$ since $B$ is strictly
increasing.

\item[(b)] If $k<i$ and $n-k\geq j$, then
$a_k'=a_k>b_{n+1-k}=b_{n-k}'$.

\item[(c)] If $k\geq i$ and $n-k<j$, then
$a_k'=a_{k+1}>b_{n-k}=b_{n-k}'$.

\item[(d)] If $k\geq i$ and $n-k\geq j$, then
$a_k'=a_{k+1}>b_{n-k+1}=b_{n-k}'$ since $M=\textbf{C}(A,B)$ is
irreducible.
\end{enumerate}
This completes the proof. \qed

The \emph{adjoint matrix} of $M$ is defined to be the matrix
$\left((-1)^{i+j}\det(M_{ji})\right)_{i,j=1}^n$, denoted $M^*$.
The \emph{rank} of a matrix $M$ is the maximum number of linearly
independent rows or columns of the matrix, denoted ${\rm rk}(M)$.
For an $n\times n$ square matrix $M$, we have the following
relationship between ${\rm rk}(M)$  and ${\rm rk}(M^*)$:
\begin{equation}\label{bd-prop}
{\rm rk}(M^{*})=\left\{
\begin{array}{ll}
n, & \mbox{ if ${\rm rk}(M)=n$,}\\[8pt]
1, & \mbox{ if ${\rm rk}(M)=n-1$},\\[8pt]
0, & \mbox{ if ${\rm rk}(M)<n-1$}.
\end{array}
\right.
\end{equation}

The restricted Cauchy matrix $M=\mathbf{C}(A,B)$ reduces to the
classical Cauchy matrix when $a_n>b_n$. In this case, the
determinant $\det(M)$ is given by the following well known formula
\begin{equation}\label{cauchy-det}
\det\left(\frac{1}{a_i-b_j}\right)_{i,j=1}^n=\prod_{i<j}(a_i-a_j)\prod_{i<j}(b_j-b_i)\prod_{i,j}\frac{1}{a_i-b_j}.
\end{equation}
Since $A$ is strictly decreasing and $B$ is strictly increasing,
the above determinant is  positive.

We now proceed to give the proof of Theorem \ref{mainthm}.

\noindent \textit{Proof of Theorem \ref{mainthm}.} We apply
induction on  $n$. The cases of $n=1,2$ are clear. Suppose that
Theorem \ref{mainthm} holds for matrices of order less than $n$.
We will prove that the theorem is also true for matrices of order
$n$.

If $M$ is a reducible Cauchy matrix, then there exists an integer
$k$ greater than or equal to 2 such that $a_k<b_{n+2-k}$. Now $M$
has the following block decomposition
$$
\begin{pmatrix}
M_1 & M_2\\
M_3 & M_4
\end{pmatrix}
,$$ where $M_1$ is a $(k-1)\times (n-k+1)$ matrix, $M_2$ is the
$(k-1)\times (k-1)$ restricted Cauchy matrix
$\textbf{C}((a_1,\ldots,a_k),(b_{n-k+1},\ldots,b_n))$, $M_3$ is
the $(n-k+1)\times (n-k+1)$ restricted Cauchy matrix
$\textbf{C}((a_{k+1},\ldots,a_n),(b_{1},\ldots,b_{n-k}))$, and
$M_4$ is an  $(n-k+1)\times (k-1)$  zero block. Thus we get
$$\det(M)=(-1)^{\omega(M_4)}\det(M_2)\det(M_3).$$ By induction
the sign of $\det(M_2)$ is $(-1)^{\omega(M_2)}$, and the sign of
$\det(M_3)$ is $(-1)^{\omega(M_3)}$. Since $\omega(M_1)=0$, the
sign of $\det(M)$ equals
$$(-1)^{\omega(M_4)+\omega(M_3)+\omega(M_2)}=(-1)^{\omega(M)}.$$

We now suppose that $M=(m_{ij})_{i,j=1}^n$ is an irreducible
Cauchy matrix. If $M$ has no zero entry, then the theorem is true
because of \eqref{cauchy-det}. If $\omega(M)>0$, then we consider
the following block decomposition of $M$
$$
\begin{pmatrix}
M_1' & M_2'\\
M_3' & 0
\end{pmatrix}
,$$ where $M_1'$ is an $(n-1)\times (n-1)$ restricted Cauchy
matrix, $M_2'$ is an $(n-1)\times 1$ column vector, $M_3'$ is a
$1\times (n-1)$ row vector. By Lemma \ref{minorcauchy}, we see
that the minors $M_{11}, M_{nn},M_{1n},M_{n1}$ are also restricted
Cauchy matrices. Consider the submatrix
$$
\begin{pmatrix}
\det(M_{11}) & (-1)^{n+1}\det(M_{n1})\\
(-1)^{n+1}\det(M_{1n}) & \det(M_{nn})
\end{pmatrix}
$$ of the adjoint matrix $M^*$. Note that the signs of $\det(M_{11}),\det(M_{n1}),\det(M_{1n})$ and $\det(M_{nn})$
are respectively
$(-1)^{\omega(M_1')+\omega(M_2')+\omega(M_3')+1},\,(-1)^{\omega(M_1')+\omega(M_2')},
\,\break(-1)^{\omega(M_1')+\omega(M_3')}$ and
$(-1)^{\omega(M_1')}$. It follows that
$$
\det\begin{pmatrix}
\det(M_{11}) & (-1)^{n+1}\det(M_{n1})\\
(-1)^{n+1}\det(M_{1n}) & \det(M_{nn})
\end{pmatrix}\neq 0.
$$ Therefore ${\rm rk}(M^*)\geq 2$. Owing to the relation
\eqref{bd-prop} between ${\rm rk}(M^*)$ and ${\rm rk}(M)$, we have
${\rm rk}(M^*)=n$, that is, $M$ is nonsingular.

It remains to show that the sign of $\det(M)$ coincides with the
number of zero entries in $M$. Without loss of generality, we may
assume that $M$ is irreducible. If $M$ does not contain any zero
entry, then $\det(M)$ is clearly positive. We now assume that $M$
contain at least one zero entry.  Note that the conditions on $A$
and $B$, for any row in $\textbf{C}(A, B)$, if there is a zero in
the $j$-column, then the entry in any column $k$ $(k>j)$ must be
zero. The same property also holds for the columns of
$\textbf{C}(A, B)$. Thus, the $(n,n)$-entry in $\textbf{C}(A, B)$
must be zero.  Since $M$ is irreducible, there exists an integer
$j:2\leq j\leq n-1$ such that $m_{nj}\neq 0$, but
$m_{n,j+1}=m_{n,j+2}=\cdots =m_{n,n}=0$. Let $\alpha=b_j$ and
$\beta=\min(a_{n-1},b_{j+1})$. Then the determinant $\det(M)$ can
be regarded as a continuous function of $a_n$ on the open interval
$(\alpha,\beta)$. Note that when $a_n$ varies in the open interval
$(\alpha,\beta)$, the restricted Cauchy matrix $M$ keeps the same
shape, which means that the positions of zero entries are fixed.
If $a_n=\eta$ for some $\eta\in(\alpha,\beta)$, denote the
corresponding matrix $M$ by $M_\eta$. When $a_n$ tends to $b_j$
from above, $m_{nj}$ tends to $+\infty$, and for $k<j$ the entry
$m_{nk}$ tends to $\frac{1}{b_j-b_k}$.

Since the minor $M_{nj}$ is a restricted Cauchy matrix of order
$n-1$, by Lemma \ref{minorcauchy}, the induction hypothesis
implies that $\det(M_{nj})\neq 0$. Therefore, the sign of
$\det(M)$ coincides with the sign of $(-1)^{n+j}\det(M_{nj})$ when
$a_n$ tends into $b_j$ from above. It follows that there exists
$\xi\in (\alpha,\beta)$ such that the sign of $\det(M_\xi)$
coincides with the sign of $(-1)^{n+j}\det(M_{nj})$. By induction,
the sign of $\det(M_{nj})$ equals $(-1)^{\omega(M_{nj})}$, thus
the sign of $\det(M_\xi)$ equals
$$(-1)^{n+j+\omega(M_{nj})}=(-1)^{(n-j)+\omega(M_{nj})}=(-1)^{\omega(M_\xi)}.$$
For any $\eta\in (\alpha,\beta)$, the sign of $\det(M_\eta)$
coincides with the sign of  $\det(M_\xi)$. Otherwise, there exists
a number $\zeta$ between $\xi$ and $\eta$ such that
$\det(M_\zeta)=0$, which is a contradiction. Since
$\omega(M_\xi)=\omega(M_\eta)$, we have completed the proof. \qed

\section{The zrank problem}

We assume that the reader is familiar with the notation and
terminology on partitions and symmetric functions  in \cite{S1}.
Given a partition $\lambda$ with decreasing components
$\lambda_1,\lambda_2,\ldots$, the rank of $\lambda$, denoted ${\rm
rank}(\lambda)$, is the number of $i$'s such that $\lambda_i\geq
i$. Clearly, ${\rm rank}(\lambda)$ counts the number of diagonal
boxes in the Young diagram of $\lambda$, where the Young diagram
is an array of squares in the plane justified from the top and
left corner with $\ell(\lambda)$ rows and $\lambda_i$ squares in
row $i$. A square $(i,j)$ in the diagram is the square in row $i$
from the top and column $j$ from the left. The content of $(i,j)$,
denoted $\tau(i,j)$, is given by  $j-i$.

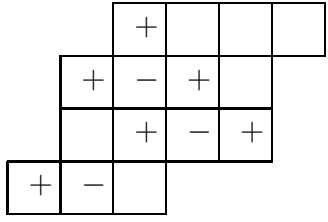
\begin{figure}[h,t]
\setlength{\unitlength}{1pt}
\begin{center}
\begin{picture}(186,100)
\put(80,85){\line(1,0){80}} \put(60,65){\line(1,0){100}}
\put(60,45){\line(1,0){80}} \put(40,25){\line(1,0){100}}
\put(40,5){\line(1,0){60}}

\put(40,5){\line(0,1){20}} \put(60,5){\line(0,1){60}}
\put(80,5){\line(0,1){80}} \put(100,5){\line(0,1){80}}
\put(120,25){\line(0,1){60}} \put(140,25){\line(0,1){60}}
\put(160,65){\line(0,1){20}}

\put(48,13){$+$}\put(68,13){$-$}\put(68,53){$+$}\put(88,33){$+$}
\put(88,53){$-$}\put(108,53){$+$}\put(88,73){$+$}
\put(108,33){$-$}\put(128,33){$+$}

\end{picture}
\end{center}
\caption{Outside and inside diagonals of
$(6,5,5,3)/(2,1,1)$}\label{in-out}
\end{figure}

Given two partitions $\lambda$ and $\mu$, if for each $i$ we have
$\lambda_i\geq \mu_i$, then the skew partition $\lambda/\mu$ is
defined to be the diagram obtained from the diagram of $\lambda$
by removing the diagram of $\mu$ at the top-left corner. A border
strip is a connected skew partition with no $2\times 2$ squares.
Nazarov and Tarasov \cite{NT} generalized rank of ordinary
partitions to skew partitions in the following way: A square
$(i,j)$ is called an \emph{inner corner} of $\lambda/\mu$, if
$(i,j),(i,j-1),(i-1,j)\in\lambda/\mu$ but $(i-1,j-1)\not\in
\lambda/\mu$; a square $(i,j)$ is called an \emph{outer corner
box} of $\lambda/\mu$, if $(i,j)\in\lambda/\mu$ but
$(i-1,j-1),(i,j-1),(i-1,j)\not\in \lambda/\mu$; the \emph{inner
diagonal} is composed of all the boxes $(i+p,j+p)\in\lambda/\mu$
if $(i,j)$ is an inner corner; the \emph{outer diagonal} is
composed of all the boxes $(i+p,j+p)\in\lambda/\mu$ if $(i,j)$ is
an outer corner; let $d^{+}$ be the number of boxes on all outer
diagonals, and let $d^{-}$ be the number of boxes on all inner
diagonals; then the \emph{rank} of $\lambda/\mu$, denoted ${\rm
rank}(\lambda/\mu)$, is the difference $d^{+}-d^{-}$. For example,
${\rm rank}((6,5,5,3)/(2,1,1))=3$, as illustrated in Figure
\ref{in-out}.

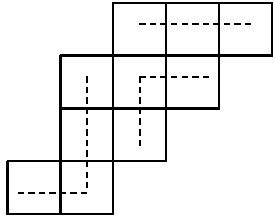
\begin{figure}[h,t]
\setlength{\unitlength}{1pt}
\begin{center}
\begin{picture}(186,100)
\put(80,85){\line(1,0){60}} \put(60,65){\line(1,0){80}}
\put(60,45){\line(1,0){60}} \put(40,25){\line(1,0){60}}
\put(40,5){\line(1,0){40}}

\put(40,5){\line(0,1){20}} \put(60,5){\line(0,1){60}}
\put(80,5){\line(0,1){80}} \put(100,25){\line(0,1){60}}
\put(120,45){\line(0,1){40}} \put(140,65){\line(0,1){20}}

\put(44,13){$\line(1,0){2}$} \put(48,13){$\line(1,0){2}$}
\put(52,13){$\line(1,0){2}$} \put(56,13){$\line(1,0){2}$}
\put(60,13){$\line(1,0){2}$} \put(64,13){$\line(1,0){2}$}
\put(68,13){$\line(1,0){2}$} \put(70,15){$\line(0,1){2}$}
\put(70,19){$\line(0,1){2}$} \put(70,23){$\line(0,1){2}$}
\put(70,27){$\line(0,1){2}$} \put(70,31){$\line(0,1){2}$}
\put(70,35){$\line(0,1){2}$} \put(70,39){$\line(0,1){2}$}
\put(70,43){$\line(0,1){2}$} \put(70,47){$\line(0,1){2}$}
\put(70,51){$\line(0,1){2}$} \put(70,55){$\line(0,1){2}$}

\put(90,31){$\line(0,1){2}$} \put(90,35){$\line(0,1){2}$}
\put(90,39){$\line(0,1){2}$} \put(90,43){$\line(0,1){2}$}
\put(90,47){$\line(0,1){2}$} \put(90,51){$\line(0,1){2}$}
\put(90,55){$\line(0,1){2}$} \put(90,57){$\line(1,0){2}$}
\put(94,57){$\line(1,0){2}$} \put(98,57){$\line(1,0){2}$}
\put(102,57){$\line(1,0){2}$} \put(106,57){$\line(1,0){2}$}
\put(110,57){$\line(1,0){2}$} \put(114,57){$\line(1,0){2}$}

\put(90,77){$\line(1,0){2}$} \put(94,77){$\line(1,0){2}$}
\put(98,77){$\line(1,0){2}$} \put(102,77){$\line(1,0){2}$}
\put(106,77){$\line(1,0){2}$} \put(110,77){$\line(1,0){2}$}
\put(114,77){$\line(1,0){2}$} \put(118,77){$\line(1,0){2}$}
\put(122,77){$\line(1,0){2}$} \put(126,77){$\line(1,0){2}$}
\put(130,77){$\line(1,0){2}$}

\end{picture}
\end{center}
\caption{A minimal border strip decomposition of
$(5,4,3,2)/(2,1,1)$}\label{min}
\end{figure}

Stanley \cite{Stanley} gave several characterizations of ${\rm
rank}(\lambda/\mu)$. The first characterization is based on the
border strip decomposition of the skew diagram. Stanley proved
that ${\rm rank}(\lambda/\mu)$ is the smallest number $r$ such
that $\lambda/\mu$ is a disjoint union of $r$ border strips. As we
see from Figure \ref{min}, \break ${\rm
rank}((5,4,3,2)/(2,1,1))=3$.

Recall that the Jacobi-Trudi identity for the skew Schur function
states that
\begin{equation}
s_{\lambda/\mu}=\det\left(h_{\lambda_i-\mu_j-i+j}\right)_{i,j=1}^{\ell(\lambda)},
\end{equation}
where $h_k$ denotes the $k$-th complete symmetric function,
$h_0=1$ and $h_k=0$ for $k<0$. Let $J_{\lambda/\mu}$ be the matrix
which appears in the Jacobi-Trudi identity. Stanley defined the
jrank of $\lambda/\mu$, denoted ${\rm jrank}(\lambda/\mu)$, by the
number of rows of $J_{\lambda/\mu}$ which do not contain $1$'s,
and proved that ${\rm jrank}(\lambda/\mu)={\rm
rank}(\lambda/\mu)$. For example,
$$
J_{(6,5,5,3)/(2,1,1)}=\begin{pmatrix} h_4 & h_6 & h_7 & h_9\\
h_2 & h_4 & h_5 & h_7\\
h_1 & h_3 & h_4 & h_6\\
0 & 1 & h_1 & h_3
\end{pmatrix}.
$$

\begin{figure}[h,t]
\setlength{\unitlength}{1pt}
\begin{center}
\begin{picture}(186,100)
\put(80,85){\line(1,0){60}} \put(60,65){\line(1,0){80}}
\put(60,45){\line(1,0){60}} \put(40,25){\line(1,0){60}}
\put(40,5){\line(1,0){40}}

\put(40,5){\line(0,1){20}} \put(60,5){\line(0,1){60}}
\put(80,5){\line(0,1){80}} \put(100,25){\line(0,1){60}}
\put(120,45){\line(0,1){40}} \put(140,65){\line(0,1){20}}

\put(142, 73){{1}}\put(130, 55){0}\put(122, 53){{1}} \put(110,
35){0}\put(102, 33){{1}} \put(90, 15){0}\put(82, 13){{1}} \put(70,
-5){0} \put(50, -5){0}

\put(75, 73){{1}}\put(130, 90){0}\put(55, 53){{1}} \put(110,
90){0}\put(55, 33){{1}} \put(90, 90){0}\put(35, 13){{1}} \put(70,
67){0} \put(50, 27){0}
\end{picture}
\end{center}
\caption{The reduced code of $(5,4,3,2)/(2,1,1)$}\label{bound}
\end{figure}
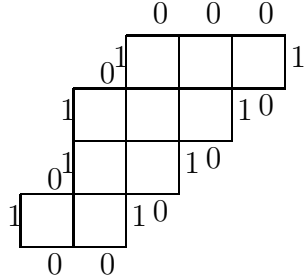

The third characterization of ${\rm rank}(\lambda/\mu)$ involves
the reduced code of $\lambda/\mu$, denoted ${\rm c}(\lambda/\mu)$.
The reduced code ${\rm c}(\lambda/\mu)$ is also known as the
\emph{partition sequence} of $\lambda/\mu$ \cite{B1, B2}. Consider
the two boundary lattice paths of the diagram of $\lambda/\mu$
with steps $(0,1)$ or $(1,0)$ from the bottom-leftmost point to
the top-rightmost point. Replacing each step $(0,1)$ by $1$ and
each step $(1,0)$ by $0$,  we obtain two binary sequences by
reading the lattice paths from the bottom-left corner to the
top-right corner. Denote the top-left binary sequence by
$g_1,\,g_2,\,\ldots,\,g_k$, and the bottom-right binary sequence
by $g_1',\,g_2',\,\ldots,\,g_k'$. The reduced code ${\rm
c}(\lambda/\mu)$ is defined by the two-row array
$$
\begin{array}{cccc}
g_1 & g_2 & \cdots & g_k\\
g_1' & g_2' & \cdots & g_k'
\end{array}.
$$
For example, the reduced code of the skew partition
$(5,4,3,2)/(2,1,1)$ in Figure \ref{bound} is
$$
\begin{array}{ccccccccc}1 & 0 & 1 & 1 & 0 &
1 & 0 & 0 & 0\\
0 & 0 & 1 & 0 & 1 & 0 & 1 & 0 & 1
\end{array}.
$$

Stanley proved that the rank of a skew partition is also equal to
the number of columns ${1 \atop 0}$ of ${\rm c}(\lambda/\mu)$, as
we see from Figure \ref{bound}.

In \cite{Stanley} Stanley introduced the notion of {\em zrank} of
a skew partition, and proposed the problem whether the zrank is
always equal to the rank for any skew partition.  Let
$s_{\lambda/\mu}(1^t)$ denote the skew Schur function
$s_{\lambda/\mu}$ evaluated at $x_1=\cdots=x_t=1,\,x_i=0$ for
$i>t$. The zrank of $\lambda/\mu$, denoted ${\rm
zrank}(\lambda/\mu)$, is the exponent of the largest power of $t$
dividing $s_{\lambda/\mu}(1^t)$.

The following equivalence was established by Yan, Yang and Zhou
\cite{YYZ}:

\begin{theo} \label{eq-des} The following two statements are equivalent:

(i) The zrank and rank are equal for any  skew partition.

(ii) Any restricted integer Cauchy matrix is nonsingular.
\end{theo}

As an immediate consequence of Theorem \ref{mainthm} and Theorem
\ref{eq-des}, we have the following conclusion.

\begin{theo} For all skew partitions $\lambda/\mu$, we have
${\rm zrank}(\lambda/\mu)={\rm rank}(\lambda/\mu)$.
\end{theo}

The above theorem allows us to give other characterizations of
rank in terms of the Giambelli type matrix, which is related to
the planar decomposition of skew diagrams introduced by Hamel and
Goulden \cite{HG}.

Let us first recall the Giambelli type determinant formulas of the
skew Schur function. A border strip decomposition of $\lambda/\mu$
is said to be an \emph{outside decomposition} if every strip in
the decomposition has an initial square on the left or bottom
perimeter of the diagram and a terminal square on the right or top
perimeter. For an outside decomposition, Chen, Yan and Yang
\cite{CYY} introduced the notion of the cutting strip and obtained
a transformation theorem on the Giambelli-type determinantal
formulas for the skew Schur function.

Suppose that $\lambda/\mu$ has $k$ diagonals, each of which is
composed of the squares having the same content. The cutting strip
of an outside decomposition is defined to be a border strip of
length $k$. Given an outside decomposition of $\lambda/\mu$, we
see that each square in the diagram can be assigned a direction in
the following way: starting with the bottom-left corner of a
strip, we say that a square of a strip has the up direction (resp.
right direction) if the next square in the strip lies on its top
(resp. to its right). Then the squares on the same diagonal of
$\lambda/\mu$ have the same direction. Based on this property, the
cutting strip $\phi$ of an outside decomposition $\mathbf{D}$ of
$\lambda/\mu$ is defined as follows: for $i=1,\,2,\,\ldots,\,k-1$
the $i$-th square in $\phi$ keeps the same direction as the $i$-th
diagonal of $\lambda/\mu$ with respect to $\mathbf{D}$. Given a
border strip $\theta$ of $\mathbf{D}$, let $p(\theta)$ denote the
lower left-hand square of $\theta$, and let $q(\theta)$ denote the
upper right-hand square. Hamel and Goulden \cite{HG} derived the
following determinantal formula.

\begin{theo}[{\cite[Theorem 3.1]{HG}}]\label{schur-dec}
For an outside decomposition $\mathbf{D}$ with $k$ border strips
$\theta_1,\,\theta_2,\,\ldots,\,\theta_k$, we have
\begin{equation}\label{hg-det}
s_{\lambda/\mu}=\det\left(s_{[\tau(p(\theta_i)),\tau(q(\theta_j))]}\right)_{i,j=1}^{k},
\end{equation}
where for any two integers $\alpha,\beta$, a strip
$[\alpha,\beta]$ is defined by the following rule: if $\alpha\leq
\beta$, then let $[\alpha,\beta]$ be the segment of $\phi$ from
the square with content $\alpha$ to the square with content
$\beta$; if $\alpha=\beta+1$, then let $[\alpha,\beta]$ be the
empty strip and $s_{[\alpha,\beta]}=1$; if $\alpha>\beta+1$, then
$[\alpha,\beta]$ is undefined and $s_{[\alpha,\beta]}=0$. The
content function $\tau$ is defined on the original skew diagram.
\end{theo}

Denote the matrix appearing in \eqref{hg-det} by
$G_{\lambda/\mu}^{\mathbf{D}}$. Given an outside decomposition
$\mathbf{D}$ of $\lambda/\mu$, let ${\rm
grank}_\mathbf{D}(\lambda/\mu)$ be the number of rows in
$G_{\lambda/\mu}^{\mathbf{D}}$ that do not contain $1$'s. Then we
have the following theorem:

\begin{theo}\label{theorem-o}
For any skew partition $\lambda/\mu$ and any outside decomposition
$\mathbf{D}$ of $\lambda/\mu$, we have ${\rm
grank}_\mathbf{D}(\lambda/\mu)={\rm rank}(\lambda/\mu)$.
\end{theo}

\proof By Theorem \ref{mainthm}, we have ${\rm
rank}(\lambda/\mu)={\rm zrank}(\lambda/\mu)$. So it suffices to
prove ${\rm grank}_\mathbf{D}(\lambda/\mu)={\rm
zrank}(\lambda/\mu)$. According to the definition of ${\rm
zrank}(\lambda/\mu)$, we need to consider the terms with the
lowest degree in the expansion of
$\det(G_{\lambda/\mu}^{\mathbf{D}}(1^t))$. Suppose that the square
with content $\tau(p(\theta_i))$ lies in the $p_i$-th row of the
cutting strip $\phi$ of $\mathbf{D}$, and the square with content
$\tau(q(\theta_j))$ lies in the $q_j$-th row. For the nonempty
border strip $[\tau(p(\theta_i)),\tau(q(\theta_j))]$, it is easy
to show that
\begin{equation}\label{j1}
(t^{-1}s_{[\tau(p(\theta_i)),\tau(q(\theta_j))]}(1^t))_{t=0}=\displaystyle\frac{(-1)^{{p_i}-{q_j}}}{\tau(q(\theta_j))+1-\tau(p(\theta_i))}.
\end{equation}
Note that for any $i\neq j$ we have
$$\tau(q(\theta_i))\neq \tau(q(\theta_j)), \quad
\tau(p(\theta_i))\neq \tau(p(\theta_j))$$ subject to the
definition of the outside decomposition $\mathbf{D}$. By removing
the rows and columns with $1$'s from
$G_{\lambda/\mu}^{\mathbf{D}}$, extracting $t$ from each row
without $1$'s, and putting $t=0$, we obtain a restricted Cauchy
matrices up to permutations  of rows and columns. From Theorem
\ref{mainthm}, we get the desired equality ${\rm
grank}_\mathbf{D}(\lambda/\mu)={\rm zrank}(\lambda/\mu)$. \qed

In fact, the above theorem can be proved in another way. Given a
border strip decomposition
$\mathbf{D}=\{\theta_1,\,\theta_2,\,\ldots,\,\theta_m\}$ of
$\lambda/\mu$, let
$$P_{\mathbf{D}}=\{\tau(p(\theta_1)),\,\tau(p(\theta_2)),\,\ldots,
\tau(p(\theta_m))\}$$ and
$$Q_{\mathbf{D}}=\{\tau(q(\theta_1))+1,\,\tau(q(\theta_2))+1,\,\ldots,
\tau(q(\theta_m))+1\}.$$ The following theorem is implicit in
\cite{YYZ}:

\begin{theo}\label{prop-o} For any border strip decomposition $\mathbf{D}$,
the two sets $P_{\mathbf{D}}- Q_{\mathbf{D}}$ and
$Q_{\mathbf{D}}-P_{\mathbf{D}}$  are independent of the border
strip decomposition $\mathbf{D}$, hence uniquely determined by the
skew shape $\lambda/\mu$.
\end{theo}

\noindent \textbf{Remark.} We omit the proof here, since it is
similar to the proof of \cite[Proposition 3.1]{YYZ}. The border
strip decomposition $\mathbf{D}$ may not be an outside
decomposition. As shown by Yan, Yang and Zhou \cite{YYZ}, these
two sets are related to the noncrossing interval sets of a given
skew partition. If $\mathbf{D}$ is a minimal border strip
decomposition, the intersection $P_{\mathbf{D}}\cap
Q_{\mathbf{D}}$ is the empty set. Otherwise, the cardinality of
the intersection $P_{\mathbf{D}}\cap Q_{\mathbf{D}}$ is equal to
the number of rows containing $1$'s in
$G_{\lambda/\mu}^{\mathbf{D}}$. From this point of view, Theorem
\ref{prop-o} is more general than Theorem \ref{theorem-o}.

\section{The factorial Cauchy matrices}

Before defining factorial Cauchy matrices and inverse binomial
matrices, let us review some background on double Schur functions.
Let  $X=\{x_1,x_2,\ldots\}$ and $Y=\{y_1,y_2,\ldots\}$ be two sets
of variables. For a positive integer $k$, we set
\begin{equation}
(x_i|Y)_k=\prod_{1 \leq j \leq k} (x_i-y_j),
\end{equation}
and define $(x_i|Y)_0=1$. Taking $y_i=i-1$, we obtain the falling
factorial $(x_i)_k  = x_i(x_i-1)\cdots(x_i-k+1).$ Taking
$y_i=1-i$, we get the rising factorial $(x_i)^k  =
x_i(x_i+1)\cdots(x_i+k-1).$

Now we give the two equivalent definitions of the double Schur
function $S_\lambda(X, Y)$. The first definition of $S_\lambda(X,
Y)$ is a determinantal form. Let $\lambda=(\lambda_1, \lambda_2,
\ldots, \lambda_n)$. We have
\begin{equation}\label{def-alg}
 S_\lambda(X, Y) = \frac{\det \biggl(
(x_i | Y)_{\lambda_j+n-j}
  \biggr)_{i,j=1}^n}{\Delta (X)},
\end{equation}
where $\Delta (X)$ is the Vandermonde determinant in $x_1, x_2,
\ldots, x_n$:
\[ \Delta(X) = \prod_{i<j} (x_i-x_j). \]

The second definition of $S_\lambda(X, Y)$ is obtained by
Macdonald \cite{Macdonald}, and Goulden and Green \cite{GG}. Chen,
Li and Louck \cite{CLL} obtained this combinatorial definition
using the lattice path methodology. Recall that a
\emph{semistandard Young tableau} $T$ of shape $\lambda$ is a
configuration of the Young diagram of $\lambda$ with positive
integers such that each row is weakly increasing and each column
is strictly increasing. Given a Young tableau $T$ and a cell
$\alpha$ of $T$, let $T({\alpha})$ be the number filled in the
cell  $\alpha$. The combinatorial definition of $S_\lambda(X, Y)$
is as follows.

\begin{theo} \label{theo-fact} Let $\lambda$ be a partition of
length $n$. Then
\[
S_\lambda(X, Y)=\sum_{T}\prod_{\alpha\in
T}(x_{T(\alpha)}-y_{T(\alpha)+\tau(\alpha)}),
\]
summing over all column strict  tableaux $T$  on $\{1, 2, \ldots,
n\}$ of shape $\lambda$.
\end{theo}

We now define the factorial Cauchy matrix and the inverse binomial
matrix. Let $A=(a_1, \ldots, a_n)$ be a strictly decreasing
sequence of real numbers, and let $B=(b_1, \ldots, b_n)$ be a
strictly increasing sequence of positive integers. Suppose that
for any $i,j$ we have $a_i> b_{n+1-i}-1$ and $a_i\neq b_j-1$. Then
we define a { matrix} $\textbf{F}(A,B)=(f_{ij})_{i,j=1}^n$ by
setting
$$
f_{ij}=\left\{
\begin{array}{ll}
{\displaystyle \frac{1}{(a_i)_{b_j}}}, & \mbox{ if $a_i> b_j-1$},\\[12pt]
0, & \mbox{ if $a_i<b_j-1$}.
\end{array}
\right.
$$
A matrix $M$ is called a \emph{factorial Cauchy matrix} if there
exist two  sequences $A$ and $B$ satisfying the above conditions
such that $M=\textbf{F}(A,B)$.

When $A$ is also a sequence of positive integers, we can define a
matrix $\textbf{R}(A,B)=(r_{ij})_{i,j=1}^n$ by setting
$$
r_{ij}=\left\{
\begin{array}{ll}
{\displaystyle {\binom{a_i}{b_j}}^{-1}}, & \mbox{ if $a_i\geq b_j$},\\[12pt]
0, & \mbox{ if $a_i<b_j$}.
\end{array}
\right.
$$
A matrix $M$ is called an \emph{inverse binomial  matrix} if there
exist two  sequences $A$ and $B$ satisfying the above conditions
such that $M=\textbf{R}(A,B)$.

Suppose that the factorial Cauchy matrix $\textbf{F}(A,B)$ has no
zero entries, i.e., the two sequences $A$ and $B$ satisfy that
$a_i>b_j-1$ for any $i,j$. In this case, we have the following
lemma:

\begin{lemm}\label{fact-pos} Let $\textbf{F}(A,B)$ be the factorial Cauchy matrix with
$a_i>b_j-1$ for any $i,j$. Then
\begin{equation}
\det(\textbf{F}(A,B))=\frac{\Delta(X)S_{\lambda}(X,Y)}{\Pi_{k=1}^n(a_k)_{b_n}},
\end{equation}
where $\lambda_j=b_n-b_j+j-n$, $x_i=a_i-b_n+1$, and $y_j=-j+1$. In
particular, we have $\det(\textbf{F}(A,B))>0$.
\end{lemm}

\proof Since $a_i>b_j$ for any $i,j$, then we have
$$\textbf{F}(A,B)=\left(\frac{1}{(a_i)_{b_j}}\right)_{i,j=1}^n.$$
Thus
$$
\begin{array}{rcl}
\det(\textbf{F}(A,B)) & = &
\displaystyle\det\left(\frac{1}{(a_i)_{b_j}}\right)_{i,j=1}^n\\[14pt]
 & = &  \displaystyle\frac{\det\left((a_i-b_n+1)^{b_n-b_j}
 \right)_{i,j=1}^n}{\Pi_{k=1}^n(a_k)_{b_n}}\\[14pt]
 & = &  \displaystyle\frac{\Delta(X)S_{\lambda}(X,Y)}{\Pi_{k=1}^n(a_k)_{b_n}},
\end{array}
$$
where the last equality follows from the algebraic definition
\eqref{def-alg} of $S_{\lambda}(X,Y)$.

Applying Theorem \ref{theo-fact}, we have
$\det(\textbf{F}(A,B))>0$. \qed

The following theorem is the main result of this section.

\begin{theo}\label{mainthm2} Any factorial Cauchy matrix
$M=\textbf{F}(A,B)$ is nonsingular. Furthermore, the determinant
$\det(M)$ is positive if $\omega(M)$ is even; or  negative if
$\omega(M)$ is odd.
\end{theo}

\proof We use induction on $n$.  Clearly,  the theorem holds when
$n=1$ or $2$. Suppose that Theorem \ref{mainthm2} holds for
matrices of order less than $n$. We proceed to prove that it is
also true for matrices of order $n$.

If $\textbf{F}(A,B)$ is a reducible factorial Cauchy matrix, then
there exists an integer $k$ greater than or equal to 2 such that
$a_k<b_{n+2-k}-1$. Now $\textbf{F}(A,B)$ has the following block
decomposition
$$
\begin{pmatrix}
F_1 & F_2\\
F_3 & F_4
\end{pmatrix}
,$$ where $F_1$ is a $(k-1)\times (n-k+1)$ matrix, $F_2$ is a
$(k-1)\times (k-1)$ factorial Cauchy matrix, $F_3$ is an
$(n-k+1)\times (n-k+1)$ factorial Cauchy matrix, and $F_4$ is an
$(n-k+1)\times (k-1)$ zero block. So we have
$$\det(\mathbf{F}(A,B))=(-1)^{\omega(F_4)}\det(F_2)\det(F_3).$$
By the assumption, the sign of $\det(F_2)$ is
$(-1)^{\omega(F_2)}$, and the sign of $\det(F_3)$ is
$(-1)^{\omega(F_3)}$. Since $\omega(F_1)=0$, the sign of
$\det(\textbf{F}(A,B))$ equals
$$(-1)^{\omega(F_4)+\omega(F_3)+\omega(F_2)}=(-1)^{\omega(\textbf{F}(A,B))}.$$

Now suppose that $\textbf{F}(A,B)=(f_{ij})_{i,j=1}^n$ is an
irreducible factorial Cauchy matrix. If $\textbf{F}(A,B)$ has no
zero entry, then the theorem is true according to Lemma
\ref{fact-pos}. If $\omega(\textbf{F}(A,B))>0$, we consider the
following block decomposition of $M$
$$
\begin{pmatrix}
F_1' & F_2'\\
F_3' & 0
\end{pmatrix}
,$$ where $F_1'$ is an $(n-1)\times (n-1)$ factorial Cauchy
matrix, $F_2'$ is an $(n-1)\times 1$ column vector, $F_3'$ is an
$1\times (n-1)$ row vector. It is easy to see that the minors
$M_{11}, M_{nn},M_{1n},M_{n1}$ of $M=\textbf{F}(A,B)$ are also
factorial Cauchy matrices. Consider the submatrix
$$
\begin{pmatrix}
\det(M_{11}) & (-1)^{n+1}\det(M_{n1})\\
(-1)^{n+1}\det(M_{1n}) & \det(M_{nn})
\end{pmatrix}
$$ of the adjoint matrix $M^*$. Note that the signs of $\det(M_{11}),\det(M_{n1}),\det(M_{1n})$
 and $\det(M_{nn})$
are given by
$(-1)^{\omega(F_1')+\omega(F_2')+\omega(F_3')+1},\,(-1)^{\omega(F_1')+\omega(F_2')},
\,\break(-1)^{\omega(F_1')+\omega(F_3')}$ and
$(-1)^{\omega(F_1')}$. Thus
$\det(M_{11})/((-1)^{n+1}\det(M_{1n}))$ and\break
$(-1)^{n+1}\det(M_{n1})/\det(M_{nn})$ must have  different signs.
Therefore, we obtain ${\rm rk}(M^*)\geq 2$. By the relation
\eqref{bd-prop}, we have ${\rm rk}(M^*)=n$, that is, $M$ is
nonsingular.

It remains to show that the sign of $\det(M)$ coincides with the
number of zero entries in $M$. Without loss of generality, we may
assume that $M$ is irreducible. If $M$ does not contain any zero
entry, then $\det(M)$ is clearly positive. We now assume that
$M=\textbf{F}(A, B)$ contain at least one zero entry. Note that
the matrix $\textbf{F}(A, B)$ has the same distribution of zeros
as the restricted Cauchy matrix.  Thus, the $(n,n)$-entry in
$\textbf{F}(A, B)$ must be zero. Since $M$ is irreducible, there
exists an integer $j:2\leq j\leq n-1$ such that $m_{nj}\neq 0$,
but $m_{n,j+1}=m_{n,j+2}=\cdots =m_{n,n}=0$. Let $\alpha=b_j-1$
and $\beta=\min(a_{n-1},b_{j+1}-1)$. Then the determinant
$\det(M)$ can be regarded as a continuous function of $a_n$ on the
open interval $(\alpha,\beta)$. Note that when $a_n$ varies in the
open interval $(\alpha,\beta)$, the factorial Cauchy matrix $M$
keeps the same shape. If $a_n=\eta$ for some
$\eta\in(\alpha,\beta)$, denote the corresponding factorial Cauchy
matrix $M$ by $M_\eta$. When $a_n$ tends to $\alpha=b_j-1$ from
above, $m_{nj}$ tends to $+\infty$, and for $k<j$ the entry
$m_{nk}$ tends to $\frac{1}{(b_j-1)_{b_k}}$.

Since the minor $M_{nj}$ is a factorial Cauchy matrix of order
$n-1$, the induction hypothesis implies that $\det(M_{nj})\neq 0$.
Therefore, the sign of $\det(M)$ coincides with the sign of
$(-1)^{n+j}\det(M_{nj})$ when $a_n$ tends to $b_j-1$ from above.
It follows that there exists $\xi\in (\alpha,\beta)$ such that the
sign of $\det(M_\xi)$ coincides with the sign of
$(-1)^{n+j}\det(M_{nj})$. By induction, the sign of $\det(M_{nj})$
equals $(-1)^{\omega(M_{nj})}$, thus the sign of $\det(M_\xi)$
equals
$$(-1)^{n+j+\omega(M_{nj})}=(-1)^{(n-j)+\omega(M_{nj})}=(-1)^{\omega(M_\xi)}.$$
For any $\eta\in (\alpha,\beta)$, the sign of $\det(M_\eta)$
coincides with the sign of  $\det(M_\xi)$. Otherwise, there exists
a number $\zeta$ between $\xi$ and $\eta$ such that
$\det(M_\zeta)=0$, which is a contradiction. Since
$\omega(M_\xi)=\omega(M_\eta)$, we have completed the proof. \qed

\begin{coro} Any inverse binomial matrix
$M=\textbf{R}(A,B)$ is nonsingular. Furthermore, the determinant
$\det(M)$ is positive if $\omega(M)$ is even; or negative if
$\omega(M)$ is odd.
\end{coro}

\proof Note that
$$
\det(\textbf{R}(A,B))=\det(r_{ij})_{i,j=1}^n=\det(\textbf{F}(A,B))\prod_{i=1}^n
b_j!.
$$
Applying Theorem \ref{mainthm2}, we get the desired assertion.
\qed

 \vspace{.2cm} \noindent{\bf Acknowledgments.} This work was supported by
 the 973 Project on Mathematical
Mechanization, the Ministry of Education, the Ministry of Science
and Technology, and the National Science Foundation of China.

\end{document}